\documentstyle[11pt]{article}
\begin{document}

\begin{center}
    {\Large {\bf Mukai Flop and Ruan Cohomology}}
\end{center}
\begin{center} Jianxun Hu \footnote[1]{supported in part by NSF of China(10171114 and 10231050)}

\end{center}

\begin{center}\small Department of Mathematics, Zhongshan University,\\
                        Guangzhou 510275
\end{center}

\begin{center}Wanchuan Zhang
\end{center}

\begin{center}\small Department of Mathematics, University of Wisconsin-Madison,\\
                     Madison WI 53706, USA
\end{center}

\begin{center}
\begin{minipage}{130mm}
\vskip 0.5cm
\begin{center}
     {\bf Abstract}
\end{center}

   Suppose that two compact manifolds $X, X'$ are connected by a sequence of Mukai flops. In this paper, 
we construct a ring  isomorphism between cohomology ring of $X$ and $X'$. Using the localization technique, we prove that the quantum corrected products on $X, X'$ are the ordinary intersection products. 
 Furthermore, $X, X'$ have isomorphic Ruan cohomology. i.e. we proved the cohomological minimal model conjecture proposed by Ruan.


\vskip 0.5cm
{\bf 2000 AMS Classification:} Primary 14N35, Secondary 53D45.

\vskip 0.5cm

\end{minipage}
\end{center}

\section{Introduction}

After the mathematical foundation of quantum cohomology was
established during last decade, see \cite{RT1}, now the focus of the research is on
its computation and application. We think that the fundamental problem in
quantum cohomology is the {\bf quantum naturality problem}\cite{R1, R2,H}:
{\it Define `` morphism'' of symplectic manifolds so that quantum cohomology
is natural.} Qin and Ruan \cite{QR} showed that the quantum cohomology is not
natural for fibrations. Their results also shows that possible ``morphism'' must
be very rigid. The existence of these rigid morphisms will set apart quantum
cohomology from ordinary cohomology and gives it its own identity. Although this
result let us feel depressed, the result of \cite{LR} discovers some amazing relations
between quantum cohomology and birational geometry. Their result said that threefolds which
are connected by a sequence of flops have isomorphic quantum cohomology. This gives us some suggestions to
look for the suitable ``morphism" from some birational transformations.

In the study of higher dimensional algebraic geometry, the famous ``minimal model program''
initiated by Mori is one of the main research topics. So far the existence problem of minimal models
is still completely open in dimensions higher than three. Moreover, in contrast to the two
dimensional case, the minimal model is not unique in higher dimensions. It is then an important question to
see what kind of invariants are shared by all the birationally equivalent minimal models, and more
generally, are preserved under certain elementary birational transformations.

It is well-known that the crepant resolutions are not unique in dimensions higher than three. But Wang \cite{W}
showed that the different crepant resolutions are connected by ``K-equivalence''. Two smooth complex manifolds
$X,Y$ are {\bf K-equivalent} if and only if there is a common resolution $\phi: Z \longrightarrow X$ and 
$\varphi : Z \longrightarrow Y$ such that $\phi^* K_X = \varphi^* K_Y$. Batyrev \cite{B} and Wang \cite{W} showed 
that two K-equivalent projective manifolds have the same betti number. It is natural to ask if they have the same 
cohomology ring structure. Unfortunately, they usually have different ring structures. About this problem,
Wang \cite{W, W1} proposed his following conjecture:

{\bf Wang's conjecture:} For K-equivalent manifolds under birational map $f: X \cdots \rightarrow X'$, there 
is a naturally attached correspondence $T \in A^{\dim X}(X \times X')$ of the form $T = \overline{\Gamma}_f + 
\sum _i T_i$ with $\overline{\Gamma}_f \subset X \times X'$ the cycle of graph closure of $f$ and with $T_i$'s 
being certain degenerate correspondences (i. e. $T_i$ has positive dimensional fibers when projecting to $X$ or
$X'$) such that $T$ is an isomorphism of Chow motives. 

In other words, Wang's conjecture implies that for K-equivalent manifolds $X$ and $X'$,  the canonical morphism 
$\varphi_* \phi^* : H^*(X, Q) \longrightarrow H^*(X', Q)$ gives rise to an isomorphism with some modification in 
the middle dimension. In the case of hyperk\"ahler manifolds, using Bishop's theorem \cite{Bi},
Huybrechts \cite{Huy1, Huy2} proved this conjecture by showing the existence of the correction 
cycles $T_i$. In this paper, for arbitrary projective manifolds connected by Mokai flops,
we proved that $\varphi_*\phi^*$ gives rise to an isomorphism of cohomology rings of $X, X'$ with an explicit expression of the correction cycles $T_i$(see the 
definition of the map $T$ in section 3).

We will concentrated our attention on a special kind of birational
transformations--{\bf Mukai Flops} \cite{Mukai}. Here we first recall the definition
of certain known flops.The simplest type of flops are called
{\bf ordinary flops}. An ordinary ${\bf P}^r$-flop(or simply ${\bf P}^r$-flop) $f: X \longrightarrow X'$ 
is a birational map such that the exceptional set $Z \subset X$ has a ${\bf P}^r$-bundle structure 
$\varphi : Z \longrightarrow S$ over some smooth variety S and the normal bundle $N_{Z/X}$ is isomorphic 
to ${\cal O}(-1)^{r+1}$ when restricting to any fiber of $\varphi$. The map $f$ and the space $X'$ are then 
obtained by first blowing up $X$ along $Z$ to get $Y$, with exceptional divisor $E$ a ${\bf P}^r \times 
{\bf P}^r$-bundle over $S$, then blowing down $E$ along another fiber direction. Ordinary ${\bf P}^r$-flops 
are also called classical flops. Three dimensional classical flops are the most well-known Atiyah flops over 
$(-1,-1)$ rational curves.

Another important example is the { \bf Mukai flops} $f: X \longrightarrow X'$. In this case it is required 
that the exceptional set $Z \subset X$ is of codimension $r$ and has a ${\bf P}^r$-bundle structure 
$\varphi : Z={\bf P}_S(F)\longrightarrow S$ (for some rank $r+1$ vector bundle $F$) over a smooth base $S$, 
moreover the normal bundle $N_{Z/X}\equiv T^*_{Z/S}$, the relative cotangent bundle of $\varphi$. To get $f$, 
one first blows up $X$ along $Z$ to get $\phi: Y \longrightarrow X$ with exceptional divisor $E=
{\bf P}_S(T^*_{Z/S})\subset {\bf P}_Z(F)\times {\bf P}_S(F^*)$ as the incidence variety. The first projection 
corresponds to $\phi$ and one may contract $E$ through the second projection to get $\phi':Y \longrightarrow X'$.

In this paper, we will only consider the following simplest Mukai flops:

{\bf Definition:} Let $X$ be a projective manifold of complex dimension $2n$. A Mukai flop from $(X,Z)$ to 
$(X', Z')$ is the following birational transformation
$$
   \begin{array}{rcl}
        & E \subset Y&\\
        &\phi \swarrow \searrow \varphi&\\
        Z \cong {\bf P}^n \subset X& \cdots \longrightarrow & X' \supset
        Z' \cong ({\bf P}^n)^*
   \end{array}
$$
where $E$ is the incidence correspondence between $Z$ and $Z'$. We also call $X$ and $X'$ are connected by 
a {\bf Mukai flop}.

Throughout this paper, we will call this simplest Mukai flops as Mukai flops. 

In the study of birational geometry, one of the most important problems is to find that what kind of cohomology is 
preserved by K-equivalene. For this purpose, Ruan \cite{R3} proposed

{\bf Quantum Minimal Model Conjecture:} Two K-equivalent projective manifolds have the same quantum cohomology.

In dimensions higher than three, quantum minimal model conjecture seems to be a difficult problem. We think the
difficulty comes from the fact we used all quantum information involving the GW-invariants. So Ruan proposed 
that we should consider another kind of cohomology with a minimal set of quantum information involving the 
GW-invariants of exceptional rational curves. We call this new cohomology as Ruan Cohomology, and  will 
give its definition in section 2. In section 4, we will also prove Ruan cohomology is invariant under Mukai flops.

Our main theorem in this paper is

{\bf Theorem:} Two compact projective manifolds which are connected by a sequence of Mukai flops have 
isomorphic cohomology and Ruan cohomology.

We will divide the proof of the theorem into two cases: ordinary cohomology and Ruan cohomology. In section 3, 
we will prove that $X, X'$ have isomorphic cohomology, see Theorem 3.2. In section 4, we will prove that for 
$X, X'$ the quantum correction all vanish. So they have isomorphic Ruan cohomology, see Theorem 4.4.

{\bf Acknowledgements:} We would like to thank Prof. Yongbin Ruan, Prof. Weiping Li for their many-hour-long
suggestive discussion and encouragement. The first author would like to Prof. Banghe Li and Yaqing Li for their discussion 
and help during my visiting The Hongkong University of Science and Technology. Thanks also to 
the organizer of the satellite conference`` Stringy orbifolds" of ICM2002 in Chengdu  for inviting us 
to announce our results.

\section{Ruan Cohomology}

In \cite{R3}, Ruan defined his quantum corrected cohomology with respect to a 
birational map. Suppose that $X, X'$ are K-equivalent and $\pi : X \cdots\rightarrow
X'$ is the birational map. Denote by $\pi ^{-1}:X' \cdots\rightarrow
X$ the inverse birational transformation of $\pi$. Let $A_1,
\cdots, A_k$ be an integral basis of the homology classes of
exceptional effective curves. We call $\pi$ nondegenerate if $A_1,
\cdots, A_k$ are linearly independent. Then the homology class of
any exceptional effective curve can be written as $A= \sum_i a_i
A_i$ for $a_i\geq 0$. For each $A_i$, we assign a formal variable
$q_i$. Then $A$ corresponds to $q_1^{a_1}\cdots q_k^{a_k}$.  We
define a 3-point function
\begin{equation}
      <\alpha, \beta, \gamma>_{qc}(q_1,\cdots, q_k) =
      \sum_{a_1,\cdots, a_k} \Psi^X_A(\alpha, \beta,
      \gamma)q_1^{a_1}\cdots q_k^{a_k}, 
\end{equation}
where $\Psi^X_A(\alpha,\beta,\gamma)$ is Gromov-Witten invariant
and $qc$ stands for the quantum correction and $\alpha, \beta, \gamma\in H^*X$. 
We view $<\alpha,\beta,\gamma>_{qc}(q_1,\cdots,q_k)$ as analytic function
of $q_1,\cdots, q_k$ and set $q_i = -1$ and let
\begin{equation}
   <\alpha, \beta,\gamma>_{qc} = <\alpha, \beta,
   \gamma>_{qc}(-1,\cdots,-1).
\end{equation}
We define a quantum corrected triple intersection
$$
  <\alpha, \beta,\gamma>_{\pi} = <\alpha,\beta,\gamma> +
  <\alpha,\beta,\gamma>_{qc},
$$
where $<\alpha, \beta,\gamma>= \int _X \alpha \cup \beta \cup \gamma$ is
the ordinary triple intersection. Then we define the quantum
corrected  product $\alpha *_{\pi}\beta$ by the equation
$$
   <\alpha *_{\pi}\beta, \gamma> = <\alpha,\beta,\gamma>_{\pi}
$$
for arbitrary $\gamma$. Another way to understand
$\alpha *_{\pi}\beta$ is as follows. Define a product as the
ordinary intersection product corrected by $\alpha *_{qc}\beta$. Namely,
\begin{equation}
       \alpha *_{\pi}\beta = \alpha \cup \beta + \alpha
       *_{qc}\beta.
\end{equation}
It is easy to see that the quantum corrected product gives rise to
a ring structure on the cohomology group of $X$, Denote this
cohomology ring as $RH_{\pi}^*(X, {\bf C})$.

{\bf Definition 2.1:} Define the quantum corrected cohomology ring
$RH^*_{\pi}(X,{\bf C})$ as {\bf Ruan cohomology} of $X$.

Ruan computed some examples of his cohomology in \cite{R3,R4}.
 About
this cohomology, Ruan \cite{R3} proposed the following conjecture

{\bf Cohomological Minimal Model Conjecture:} Suppose that $\pi: X\longrightarrow X'$
and its inverse $\pi^{-1}$ are nondegenerate. Then $RH^*_{\pi}(X, {\bf C})$ is isomorphic 
to $RH^*_{\pi^{-1}}(X', {\bf C})$.

{\bf Example 2.2:} The first example is the flop in dimension three. This case has been worked out
in great detail by Li-Ruan\cite{LR}. For example, they proved a theorem that quantum cohomology rings are
isomorphic under the change of the variable $q\longrightarrow \frac{1}{q}$. Notes that if we set
$q=-1$, $\frac{1}{q}= -1$. We set other quantum variables zero. Then, the quantum product becomes the 
quantum corrected product $\alpha \cup _{\pi}\beta$. Hence, Cohomological Minimal Model conjecture follows 
from LI-Ruan's theorem.  In fact, it is easy to calculate the quantum corrected product in this case and 
verify the Cohomological Minimal Model conjecture without using Li-Ruan's theorem.  

\section{Isomorphism of ordinary cohomology }

In this section, we will consider the cohomology of compact projective manifolds of complex dimension $2n$ connected by Mukai flops. Suppose that $X$ and $X'$ are compact projective manifolds of complex dimension $2n$, and
$(X,{\bf P}^n)$ and $(X', ({\bf P}^n)^*)$ are connected by a Mukai flop. Now the normal bundle 
of ${\bf P}^n$ in $X$ is its cotangent bundle $T^*{\bf P}^n$. So we have the 
following Mukai transformation
$$
      \begin{array}{rcl}
        & E \subset \tilde{X}&\\
        &\phi \swarrow \searrow \phi'&\\
        Z \cong {\bf P}^n \subset X& -- \longrightarrow & X' \supset
        ({\bf P}^n)^* \cong Z'
   \end{array}
$$
where $\tilde{X}$ is the blowup of $X$ along $Z={\bf P}^n$ and $E$ is the incidence 
correspondence between $Z$ and $Z'$, i.e. 
\begin{equation}
  \begin{array}{rcl}
      E & = & \{(P,L)\mid P \in L \}\subset {\bf P}^n \times ({\bf P}^n)^*\\
	    p \swarrow& & \searrow q \\
		P \in {\bf P}^n & & L \in ({\bf P}^n)^*.
	\end{array}
\end{equation} 

Before we prove our theorem, we want to first introduce some notations and preliminary results.
Let $X$ be a regularly embedded subscheme of a scheme $Y$ of codimension $d$ with normal bundle $N$. 
Let $A_k(X)$ be the group of $k$-cycles modulo rational equivalence on $X$. Denote by $s(X,Y)\in A_*(X)$ 
the Segre class of $X$ in $Y$, for its definition see Section 4.2 of \cite{F}, so $s(X,Y)$ is the cap 
product of the total inverse Chern class of the normal bundle with $[X]$. Let $\tilde{Y}$ be the blowup of 
$Y$ along $X$, and let $\tilde{X} = {\bf P}(N)$ be the exceptional divisor. We have a fiber square
\begin{equation}
\begin{array}{ccc}
    \tilde{X}& \stackrel{j}{\longrightarrow}& \tilde{Y}\\
	g \downarrow & & \downarrow f\\
	 X &\stackrel{\longrightarrow}{i}& Y .
	 \end{array}
\end{equation}
Since $N_{\tilde{X}}\tilde{Y} = {\cal O}(-1)$, the excess normal bundle $\xi$
is the universal quotient bundle on ${\bf P}(N)$: 
$$
  \xi = g^*N/N_{\tilde{X}}\tilde{Y} = g^*N/{\cal O}(-1).
$$
Then we have the following {\bf Blowup formula}, which is the {\bf Theorem 6.7}, see P. 116, of \cite{F}, 

 {\bf Proposition 3.1:} Let $V$ be a $k$-dimensional subvariety of $Y$, and let $\tilde{V}\subset \tilde{Y}$
 be the proper transform of $V$, i. e. the blow-up of $V$ along $V \cap X$. Then
 $$
  f^*[V] = [\tilde{V}] + j_*\{c(\xi)\cap g^*s(V \cap X, V)\}_k
 $$
 in $A_k \tilde{Y}$. In particular, for all $x \in A_k X$, 
 $$
  f^*i_*(x) = j_*(c_{d-1}(E)\cap g^*x).
 $$
 
In our proof, we will use Borel-Moore homology as a tool. Therefore we first want to briefly introduce some basics of Borel-Moore homology, see \cite{CG,F}. Borel-Moore homology can be defined using singular cohomology. If a space $X$ is imbedded as a closed subspace of ${\bf R}^n$, then we define the Borel-Moore homology with rational coefficients
$$
     H^{BM}_iX := H^{n-i}({\bf R}^n, {\bf R}^n-X)
$$
where the group on the right is relative singular cohomology with rational coefficients. From the difinition, it is easy to know if $X$ is compact then the ordinary homology of $X$ and the Borel-Moore homology of $X$ coincide. In this paper, we will reserve the symbol $H_*$ for the ordinary homology.

If $X$ is the complement of $U$ in $Y$, $i : X\longrightarrow Y$ the closed imbedding, there is a long exact sequence
\begin{equation}
  \cdots \rightarrow H^{BM}_{i+1}U\rightarrow H^{BM}_i X \stackrel{i_*}{\rightarrow}H^{BM}_i Y \stackrel{j^*}{\rightarrow} H^{BM}_i U \rightarrow H^{BM}_{i-1}X \rightarrow \cdots.
\end{equation} 

In this section, we will prove the following theorem
 
 {\bf Theorem 3.2:} Suppose that non-singular compact projective manifolds $X$ and $X'$ of complex dimension $2n$
 are connected by a sequence of Mukai flops. Then $X$ and $X' $ have isomorphic cohomology rings.
 
 {\bf  Proof:} By the Poincare duality, it is sufficient to prove that $X$ and $X'$ have isomorphic intersection rings. 
 In fact, we want to prove  the following morphism $T : H_*X \longrightarrow H_*X'$ given by
 $$
       T( \alpha):= \left\{ 
  \begin{array}{ll}
      \phi'_*\phi^* \alpha,& \mbox{if} \dim \alpha \not= 2n\\
	   \phi'_* (\phi^*\alpha + (-1)^{n+1}\alpha({\bf P}^n)[p^{-1}({\bf P}^1)]), & \mbox{if} \dim \alpha = 2n 
  \end{array}\right.
 $$
 is a ring isomorphism, where $\alpha ({\bf P}^n)$ is the topological intersection number of $\alpha$ with 
 ${\bf P}^n$ and ${\bf P}^1$ is a line in ${\bf P}^n$.  It is obvious that $T$ is a linear map.
 
 First of all, we want to prove that the restriction of $T$ to $i_* H_k ({\bf P}^n)$ is an isomorphism from 
 $i_* H_k ({\bf P}^n)$ to $i'_* H_k (({\bf P}^n)^*)$. By the linearity of $T$, we only need to prove that $T$ maps 
 a basis of $i_* H_* ({\bf P}^n)$ to a basis of $i'_* H_* (({\bf P}^n)^*)$. Since all elements in $i_* H_* ({\bf P}^n)$ are algebraic, so we may apply proposition 3.1. In our case, we have the following 
 blowup fiber square 
\begin{equation}
\begin{array}{ccc}
    E& \stackrel{j}{\longrightarrow}& \tilde{X}\\
	p \downarrow & & \downarrow \phi\\
	 {\bf P}^n &\stackrel{\longrightarrow}{i}& X.
	 \end{array}
\end{equation} 
where $i$ embedded ${\bf P}^n$ into $X$ with its cotangent bundle $N_{{\bf P}^n|X}\cong T^*{\bf P}^n$ as the 
normal bundle and $E$ is the exceptional divisor. The excess normal bundle $Q$ is the universal 
quotient bundle on $E$
$$
Q = \frac{p^*T^*{\bf P}^n}{{\cal O}_E(-1)}
$$
i. e. we have the exact sequence
$$
0 \longrightarrow {\cal O}_E(-1)\longrightarrow p^*T^*{\bf P}^n \longrightarrow Q \longrightarrow 0
$$ 
According to Proposition 3.1, we need to compute the Chern class $c_{n-1}(Q)$. Since $c(p^*T^*{\bf P}^n) = 
c(Q)c({\cal O}_E(-1))$, so we have 
\begin{eqnarray*}
c(Q) & = & \frac{c(p^*T^*{\bf P}^n)}{c({\cal O}_E(-1))}\\
    & = & \sum _{k=0}^{2n-1}\sum _{i+j = k}(-1)^i
	\left(\begin{array}{c}
	   n+1\\
	   i
	   \end{array}\right)(p^*H)^i c_1({\cal O}_E(1))^j
\end{eqnarray*}
where $H$ is the hyperplane class of ${\bf P}^n$. Therefore
\begin{eqnarray*}
 c_{n-1}(Q) & = & \sum_{i+j=n-1} (-1)^i \left (\begin{array}{c}
	   n+1\\
	   i
	   \end{array}\right )(p^*H)^ic_1({\cal O}_E(1))^j\\
	   & = & \sum _{i=0}^{n-1}\sum _{j=0}^{n-i-1}(-1)^i \left(\begin{array}{c}
	   n+1\\
	   i
	   \end{array}\right)\left(\begin{array}{c}
	   n-i-1\\
	   j
	   \end{array}\right)(q^*H^*)^{n-i-j-1}(p^*H)^{i+j}.
\end{eqnarray*}
where $H^*$ is the hyperplane class of $({\bf P}^n)^*$ and we used that $c_1({\cal O}_E(1)) = p^*H + q^*H^*$.

 Choose $i_*[{\bf P}^k]$, $k=0,\cdots,n$ as a basis of $i_* H_* ({\bf P}^n)$. For arbitrary $1 \leq k < n$,
 i. e. $x = i_*[{\bf P}^k]\in i_*H_*({\bf P}^n)$, by Proposition $3.1$,  we have
 \begin{eqnarray*}
    \phi^*(i_*[{\bf P}^k])& = & j_*\{ \sum _{i=0}^{n-1}\sum _{j=0}^{n-i-1}(-1)^i \left(\begin{array}{c}
	   n+1\\
	   i
	   \end{array}\right)\left(\begin{array}{c}
	   n-i-1\\
	   j
	   \end{array}\right)(q^*H^*)^{n-i-j-1}(p^*H)^{i+j}\cap p^*[{\bf P}^k]\}\\
	   & = & j_*\{ \sum _{i=0}^{n-1}\sum _{j=0}^{n-i-1}(-1)^i \left(\begin{array}{c}
	   n+1\\
	   i
	   \end{array}\right)\left(\begin{array}{c}
	   n-i-1\\
	   j
	   \end{array}\right)(q^*H^*)^{n-i-j-1}\cap p^*(H^{i+j}\cap [{\bf P}^k])\}.
\end{eqnarray*}
Therefore, we have
\begin{eqnarray}
\phi'_*\phi^*(i_*[{\bf P}^k]) & = & \{\sum _{i=0}^k(-1)^i \left(\begin{array}{c}
         n+1\\
		 i \end{array}\right)\left(\begin{array}{c}
		 n-i\\
		 k-i \end{array}\right)\}i'_*([({\bf P}^k)^*]) \nonumber\\
	& = & (-1)^k i'_*([({\bf P}^k)^*])
\end{eqnarray}
where we used the facts that for any $k \geq 2$ the maps $\phi' : p^*[{\bf P}^k]\longrightarrow 
({\bf P}^n)^*$ have positive dimensional fibers.

For the case $k=n$, i. e. $x= i_*({\bf P}^n)\in i_*H_*({\bf P}^n)$, we have
\begin{eqnarray*}
\phi^*(i_*[{\bf P}^n]) & = & j_*\{ \sum _{i=0}^{n-1}\sum _{j=0}^{n-i-1}(-1)^i \left(\begin{array}{c}
	   n+1\\
	   i
	   \end{array}\right)\left(\begin{array}{c}
	   n-i-1\\
	   j
	   \end{array}\right)(q^*H^*)^{n-i-j-1}(p^*H)^{i+j}\cap p^*[{\bf P}^n]\}\\
	& = & \sum_{i=0}^{n-1}(-1)^i \left(\begin{array}{c}
	  n+1\\
	  i \end{array}	\right)p^*[{\bf P}^1]\\
	& = & (-1)^{n+1}n p^*[{\bf P}^1].
\end{eqnarray*}
Therefore,
\begin{equation}
 \phi'_*\phi^*(i_*({\bf P}^n)) = (-1)^{n+1}n\phi'_*(p^*[{\bf P}^1]) = (-1)^{n+1}n i'_*({\bf P}^n)^*.
\end{equation}

Furthermore, by the definition of the map $T$, we have
\begin{equation}
    T(i_*({\bf P}^n)) = (-1)^n i'_*({\bf P}^n)^*.
\end{equation}

Since $i'_*({\bf P}^k)^*$, $k=0,\cdots, n$ is a basis of $i'_*H_*(({\bf P}^n)^*)$, so the restriction of $T $
to $i_*H_*({\bf P}^n)$ is an isomorphism from $i_*H_*({\bf P}^n)$ to  $i'_*H_*(({\bf P}^n)^*)$.

Next we want to prove that $T$ is an isomorphism of additive homology.

Denote $U:=X-{\bf P }^n$ and $U' := X'-({\bf P}^n)^*$.  Since $H^{BM}_i{\bf P}^n$ has at most one generator for all $i$, then, from $(6)$, we have the following exact sequences:
\begin{eqnarray}
   0 \longrightarrow i_*H^{BM}_k{\bf P}^n \stackrel{\subset}{\longrightarrow} H^{BM}_k X \stackrel{j^*}{\longrightarrow} 
      H^{BM}_k U \longrightarrow 0 \\
   0 \longrightarrow i_*H^{BM}_k ({\bf P}^n)^* \stackrel{\subset}{\longrightarrow} H^{BM}_k X' \stackrel{j^*}{\longrightarrow} 
      H^{BM}_k U' \longrightarrow 0. 
\end{eqnarray}
 Since $H^{BM}_k {\bf P}^n$, $H^{BM}_k X$, $H^{BM}_k U$, $H^{BM}_k ({\bf P}^n)^*$, $H^{BM}_k X'$, $H^{BM}_k U'$  all are free Abelian groups,
so we have 
\begin{eqnarray}
      H^{BM}_k X & \cong & i_*H^{BM}_k {\bf P}^n \oplus H^{BM}_k U \\
      H^{BM}_k X'& \cong & i_*H^{BM}_k({\bf P}^n)^* \oplus H^{BM}_k U'.
\end{eqnarray}
Here we used the following elementary fact from extension theory:

{\bf Proposition 3.3: (\cite{Bott}, P. 168) } In a short exact sequence of Abelian groups
$$
0 \longrightarrow A \longrightarrow B \longrightarrow C \longrightarrow 0,
$$ 	   	   
if $A$ and $C$ are free, then $ B \cong A \oplus C$.

 In fact, the previous proof shows that the restriction of $T$ to $i_*H^{BM}_k({\bf P}^n)$ is an isomorphism from $i_*H^{BM}_k({\bf P}^n)$ to $i'_*H^{BM}_k(({\bf P}^n)^*)$. On the other hand, since $\phi$ and $\phi'$
are the identity map outside ${\bf P}^n$ and $({\bf P}^n)^*$ respectively, i. e. $U \cong U'$,
the restriction of $T$ to $H^{BM}_k U$ is also an isomorphism from $H^{BM}_k U$ to $H^{BM}_k U'$. By the 
linearity of  $T$, from $(13)$ and $(14)$, we have that $T$ is an isomorphism from $H^{BM}_k X$ to $H^{BM}_k X'$ as additive groups. Since $X$ and $X'$ are compact, Therefore, $T$ also gives an isomorphism
from the ordinary homology $H_k X$ to $H_k X'$ as additive groups. 

 Now it remains to prove that $T$ preserves the multiplication, i. e. for any classes $\alpha, \beta \in H_*X$,  
 we have  
 \begin{equation}
     T(\alpha \cdot \beta) = T(\alpha)\cdot T(\beta).
\end{equation}
 
By the transverality theorem, for any homology classes $\alpha, \beta$, we may choose their representatives $M$(for $\alpha$) and $N$(for $\beta$) respectively such that they transversally intersect, i. e. $\dim (M\cap N) = \dim M + \dim N - 4n$. In the rest proof of this theorem, we will use the same symbol to denote the homology class and its representatives.
  
Since $T$ is linear and the intersection product is distributive, we only need to prove $(15)$ holds for generator
 classes. From the fact that the intersection product is a map from $H_k X \otimes H_l X$ to $H_{k+l-4n}X$, we know that
 $(15)$ holds if $\dim \alpha + \dim \beta < 4n$. Therefore, we may assume that $\dim \beta \geq 2n$. Since $U:= X-Z$ is
 isomorphic to $U':=X'-Z'$, we have that the map $T$ is the identity map on $ H_*(X-Z)$. Therefore, If at least
 one of the supports of $\alpha, \beta$ does not intersect with ${\bf P}^n$, then $(15)$ holds. Therefore, we only need to 
 consider the following four cases.
 
 {\bf Case I:} $\dim \alpha <2n$, $\beta $ is an arbitrary class.
 
 In this case, we may choose a representative submanifold $\alpha $ with support away from ${\bf P}^n$. Therefore, by the construction of the intersection product and the fact that 
 $T $ is an identity map from $H_*(U)$ to $H_*(U')$, we have
 $$
  T(\alpha\cdot \beta) = T(\alpha)\cdot T(\beta).
 $$ 
 
 {\bf Case II:} $\dim \alpha = 2n$ and $\dim \beta = 2n$.

From $(13)$ and the distributivity of intersection product, we only need to consider the case: $\alpha = i_*({\bf P}^n)$ and $\beta = i_*({\bf P}^n)$. In this case, we have  
 $$
 T(\alpha\cdot\beta) = T(-(n+1)[pt]) = -(n+1)[pt] = i'_*({\bf P}^n)^*\cdot i'_*({\bf P}^n)^* = T(\alpha)\cdot T(\beta).
 $$

 {\bf Case III:} $\dim \alpha >2n $ and $\dim \beta >2n$.
 
 Here we first prove the following claim: 
 
 {\bf Claim:}If $\phi : \tilde{X}\longrightarrow X$ is the blowup of $X$ along a 
 subvariety, then $\phi^*\alpha\cdot\phi^*\beta = \phi^*(\alpha\cdot \beta)$ for any classes 
$\alpha, \beta \in H_*X$.

In fact, by definition, we have 
\begin{eqnarray*}
    \phi^*(\alpha \cdot \beta) & = & PD \phi^* PD (\alpha \cdot \beta)
                            = PD \phi^*(PD(\alpha)\cup PD(\beta))\\
  & = & PD(\phi^*PD(\alpha)\cup \phi^*PD(\beta) = PD\phi^*PD(\alpha)\cdot PD\phi^* PD(\beta)\\
  & = & \phi^* \alpha \cdot \phi^*\beta
\end{eqnarray*}
where $PD$ stands for Poincare dual.  
 
 Since $\phi': \tilde{X'}\longrightarrow X'$ is the projection of blowup, so we have $\phi'_*\phi'^*\alpha = \alpha$ for any $\alpha \in H_*X'$. From the definition of $T$,
we have $\phi'^* T(\alpha) = \phi^*\alpha + \xi$, $\phi'^* T(\beta) = \phi^*\beta +\eta$ where $\phi'_*\xi = \phi'_*\eta = 0$, i. e. $\phi'\mid _{\xi}$ and $\phi'\mid _{\eta}$ have positive dimensional fiber. Therefore, if $\dim (\alpha \cdot \beta)\not= 2n$, from the above claim and the projection formula, we have
\begin{eqnarray*}
   T(\alpha\cdot \beta) & = & \phi'_*\phi^*(\alpha\cdot\beta)
                         =  \phi'_*\{\phi^*\alpha\cdot\phi^*\beta\}\\
	& = & \phi'_*\{\phi'^*T(\alpha)\cdot \phi'^*T(\beta)-\phi'^*T(\alpha)\cdot\eta - \phi'^*T(\beta)\cdot                  \xi + \xi\cdot\eta\}\\
	& = & T(\alpha)\cdot T(\beta).
\end{eqnarray*}

If $\dim (\alpha \cdot \beta )= 2n$, i. e., $\dim \alpha + \dim \beta = 5n$, without loss of
generality, we may assume that $\alpha \cdot \beta = k i_*{\bf P}^n$ and $\dim \beta <4n$. By
the definition of $T$ and the intersection product, we also may assume that $T(\alpha) \cdot
T(\beta) = m i'_*({\bf P}^n)^*$. Choose a $l$-dimensional class $\gamma$ where $l$ satisfies $\dim \beta + l
- 4n < 2n$ and $l < 2n$. Then from the associativity of the intersection product and Case I, we 
have the triple intersection equality.
$$
 T(\alpha\cdot\beta\cdot\gamma) = T(\alpha)\cdot T(\beta\cdot\gamma)= T(\alpha)\cdot T(\beta)\cdot T(\gamma).
$$
Since $T(\alpha \cdot \beta \cdot \gamma)= T((\alpha\cdot\beta)\cdot\gamma)= T(\alpha\cdot\beta)
\cdot T(\gamma) = (-1)^nk i'_*({\bf P}^n)^*\cdot T(\gamma)$ and $T(\alpha)\cdot T(\beta)\cdot 
T(\gamma) = m i'_*({\bf P}^n)^* \cdot T(\gamma)$, so we have $m = (-1)^n k$. Therefore $(15)$
holds.
 
{\bf Case IV:} $\alpha = i_*{\bf P}^n$, $\dim \beta >2n$ and $\beta$ transverally intersects with ${\bf P}^n$.
 
Since all odd-dimensional classes in ${\bf P}^n$ are homologous to zero, without loss of generality, we may assume that $\dim \beta$ is even. Suppose that $\gamma$ is any $(6n-\dim \beta)$-dimensional class in $H_*X$. Then the intersection product $\beta\cdot\gamma$ is a $2n$-dimensional class in $H_{2n}X$. From the associativity of the intersection product and Case II and III, 
we have the triple intersection equality
\begin{equation}
     T(\alpha\cdot\beta\cdot\gamma) = T(\alpha)\cdot T(\beta\cdot\gamma)= T(\alpha)\cdot T(\beta)\cdot T(\gamma).
\end{equation} 

Suppose that ${\bf P}^n \cdot\beta = m i_*[{\bf P}^{\frac{\dim \beta}{2} -n}]$ and $({\bf P}^n)^*\cdot T(\beta) = 
k i'_*([{\bf P}^{\frac{\dim \beta}{2} - n}]^*)$. Then by Case I we have 
\begin{eqnarray*}
   T(\alpha\cdot\beta\cdot\gamma)& = & m T(i_*[{\bf P}^{\frac{\dim \beta}{2} -n}] \cdot\gamma ) \\
   & = & m T(i_*[{\bf P}^{\frac{\dim \beta}{2} -n}])
   \cdot T(\gamma) = (-1)^{\frac{\dim \beta}{2} -n}m i'_*([{\bf P}^{\frac{\dim \beta} - n}])^*\cdot T(\gamma).
\end{eqnarray*}
On the other hand, 
$$
 T(\alpha)\cdot T(\beta)\cdot T(\gamma) = k i'_*([{\bf P}^{\frac{\dim \beta}{2} - n}])^*\cdot T(\gamma).
 $$
 Therefore we have $m=(-1)^{\frac{\dim \beta}{2} -n}k$. Therefore
 \begin{eqnarray*}
  T(\alpha\cdot\beta) & = & m T(i_*[{\bf P}^{\frac{\dim \beta}{2} -n}])= (-1)^{\frac{\dim \beta}{2} -n}k T( i_*[{\bf P}^{\frac{\dim \beta}{2} -n}])\\ 
	& = & k i'_*([{\bf P}^{\frac{\dim \beta}{2} - n}]^*) = ({\bf P}^n)^*\cdot T(\beta) = T(\alpha)\cdot T(\beta).
\end{eqnarray*}	   	
So we proved the equality $(15)$. This proves Theorem $3.2$. 

\section{Isomorphism of Ruan Cohomology}

In this section, we will study Ruan cohomologies of $X$ and $X'$. From the previous section, we know that 
in order to prove isomorphisim of Ruan cohomology for the pair ,$X$ and $X'$,  we need to calculate the quantum 
corrected product coming from exceptional effective curves on $X$ and $X'$ respectively. In fact, we will prove vanishing of the exceptional Gromov-Witten 
invariants appearing in the definition of quantum corrected product by localization technique.

\subsection{Introduction to Localization}
The calculation of the exceptional quantum product is local in nature, i.e. only a neighborhood of the embeded ${\bf P}^n$
in $X$ or $X'$ is relevant to the quantumn product with base homology being exceptional curves living in the embeded ${\bf P}^n$.
Similar local invariants appeared in the study of local mirror symmetry.  
 As explained in
\cite{CKYZ}, local mirror symmetry refers to a specialization of mirror symmetry technique to study geometry of Fano 
surfaces inside Calabi-Yau manifolds. 

Following \cite{CKYZ}, we first briefly describe the calculation setup. Let 
$\overline{\cal M}_{0,0}({\bf P},d)$ be Kontsevich's moduli space of stable maps of genus 0(could be of higher genus)
with no marked points. Denote a point in the space by $(C,f)$, where $f : C \longrightarrow {\bf P}$ (${\bf P}$ is 
some toric variety ), and $[f(C)] = d \in H_2({\bf P})$. Let $\overline{\cal M}_{0,1}({\bf P},d)$ be the same  but with one
marked point. Consider the following diagram
$$
  \overline{\cal M}_{0,0}({\bf P},d)\longleftarrow \overline{\cal M}_{0,1}({\bf P},d)\longrightarrow {\bf P},
$$
where the first arrow denotes the forgetting map $\rho : \overline{\cal M}_{0,1}({\bf P},d)\longrightarrow 
\overline{\cal M}_{0,0}({\bf P},d)$ which forgets the marked point following stablization of the 
domain curve and the second arrow denotes the evaluation map 
$ev : \overline{\cal M}_{0,1}({\bf P}, d)\longrightarrow {\bf P}$ sending $(C,f, x_1)$ to $f(x_1)$.

Let {\bf Q} be Calabi-Yau defined as the zero section of a convex bundle $V $ over ${\bf P}$ (here convex means
$H^1(C,f^*V)=0$ for any stable map $(C,f)$). Then $U_d$ is the bundle over $\overline{\cal M}_{0,0}({\bf P},d)$
defined by
$$
   U_d := \rho_*ev^*(V).
$$
The fiber of $U_d$ over a point $(C,f)$ is $H^0(C,f^*V)$. And the Kontsevich numbers (Gromov-Witten type invariant)
are defined to be
$$
   K_d := \int_{\overline{\cal M}_{0,0}({\bf P},d)}c(U_d)
$$
where $c$ is the appropriate Chern class in the context.

In case the bundle $V$ is also concave (meaning  $H^0(C,f^*V) = 0$ for any stable map $(C,f)$), there is also an 
induced bundle over the moduli space of maps whose fiber over a point $(C,f)$ is given by $H^1(C, f^*V)$. In particular
if $V$ is the normal bundle of ${\bf P}$ with respect to certain embedding of ${\bf P}$, the induced bundle is usually
called the obstruction bundle.

In the same spirit of the above setup, there is another well known  example (the multiple cover contribution) which we now 
describe.

Let $C_0 = {\bf P}^1$ be a smooth ${\bf P}^1$ embedded in a Calabi-Yau 3-fold $M$ with balanced normal bundle 
${\cal O}(-1)\oplus {\cal O}(-1)$. The moduli space of stable maps $\overline{\cal M}_{0,0}(M, d[C_0])$ has a
connected component $\overline{\cal M}_{C_0}$ isomorphic to $\overline{\cal M}({\bf P}^1, d[{\bf P}^1])$ consisting
stable d-fold covers of $C_0$. This component has  dimension $2 d -2$ while the virtual dimension is $0$. So to correctly
count the number of maps (or to define the corresponding Gromov-Witten invariant), we have to consider the obstruction
bundle $U_d$ whose fiber over $(C,f)$ is given by $H^1(C,f^*N_{C_0|M}) = {\bf C}^2 \otimes H^1(C, f^*{\cal O}(-1))$. Note that
the rank of the obstruction bundle is also $2d-2$. And the contribution of $\overline{\cal M}_{C_0}$ is given by
$$
  M_d := \int_{\overline{\cal M}({\bf P}^1, d)}c_{2d-2}(U_d).
$$   

The above definition is proposed by Kontsevich who also derived a graph summation formula for it. And the value is 
checked by Y. Manin to be $\frac{1}{d^3}$. ( there is difficulty in summing up all the contributions from admmisible
graphs). 

The essence in both examples described above is to determine and evaluate certain cohomology class (over the space of 
stable maps) which come from bundles induced from bundles over the target space. And solutions to both problems come out
of application of localization techniques. Since the target space is toric, the moduli space of maps together with the 
induced bundles inherit torus action ( action on space of maps by translating maps). Hence the classes under consideration
can be localized to the fixed points loci and become much more accessible.

In \cite{CKYZ}, the authors considered the cases where the bundle $V$ is a direct sum of line bundles, while in this 
paper we will consider the case where the target space is ${\bf P}^n$ and the bundle $V$ is the cotangent bundle of
${\bf P}^n$  which is a natural example of concave bundles. It is of interest also because it demonstrate rather different
phenomena from the examples described above. We will describe obstruction bundle induced from cotangent bundle of ${\bf P}^n$
 and define related Gromov-Witten type invariants. Surprising we will see that all these invariants are $0$.

The essential fact used in the proof is the following {\bf observation:} let $C$ be a smooth ${\bf P}^1$ mapping  onto
a line (${\bf P}^1$) inside ${\bf P}^n$ with degree d. Denote the map by $f$. Standard torus action (diagonal action)
on ${\bf P}^n$ naturally lifted to $T^*{\bf P}^n$ induces an action on the vector space $H^1(C,f^*T^*{\bf P}^n)$. 
Calculate the weights of the action, we see that there is a $0$ weight piece.

This observation of the $0$ weight piece also leads to other interesting applications. For instance, by utilizing
it, we can calculate all the Gromov-Witten invariant, hence determine the quantum cohomology ring structure of the 
projective bundle ${\bf P}(T^*{\bf P}^2\oplus {\cal O})$ over ${\bf P}^2$. Again the difficulty lies in how to sum up,
granted with the graph summation machinery developed by Kontsevich. And the simple observation we have will greatly 
simplify the summation procedure.

The rest of this section is organized as follows: In subsection $4.2$, we define our invariant and state the vanishing 
theorem. In subsection $4.3$, we introduce the Bott's residue formula and Kontsevich's graph summation formula for computing
the invariants. In subsection $4.4$, we prove our vanishing theorem and  our result about isomorphism of 
Ruan cohomology.

\subsection{Definition of invariants}

In this subsection we define our invariants. Let $\overline{\cal M}_{g,k}({\bf P}^n,d)$ be the moduli space of stable 
maps from genus $g$  curves with $k$ marked points into ${\bf P}^n$ which carries the fundamental class 
$d[{\bf P}^1]\in H_2({\bf P}^n)$. Denote a typical element in $\overline{\cal M}_{g,k}({\bf P}^n,d)$ by 
$(C,f,x_1,\cdots,x_k)$. The cotangent bundle of ${\bf P}^n$ induces an obstruction bundle over
$\overline{\cal M}_{g,k}({\bf P}^n,d)$ whose fiber at $(C,f, x_1,\cdots,x_k)$ is $H^1(C,fT^*{\bf P}^n)$. Its Euler
class ( denoted by $\Phi$) plays an important role in defining our invariants.

There are also other cohomology classes on  $\overline{\cal M}_{g,k}({\bf P}^n,d)$. For instance there is the evaluation 
maps $ev_i : \overline{\cal M}_{g,k}({\bf P}^n,d) \longrightarrow {\bf P}^n$, sending $(C,f,x_1,\cdots, x_k)$ to 
$f(x_i)$, So we can pull back cohomology classes from ${\bf P}^n$ via the evaluation maps. Also there is the forgetting 
map $\overline{\cal M}_{g,k}({\bf P}^n,d) \longrightarrow \overline{\cal M}_{g,k}$ by forgetting the map $f$ of
$(C,f,x_1,\cdots, x_k)$ where $\overline{\cal M}_{g,k}$ is the Deligne-Munford space of stable curves with $k$ marked
points. So we can also pull back classes from $\overline{\cal M}_{g,k}$.

Integrating polynomials in these classes over the moduli space $\overline{\cal M}_{g,k}({\bf P}^n,d)$, we get numbers.

In particular, if ${\bf P}^n$ is embedded in a variety $X$ with normal bundle naturally isomorphic to its cotangent bundle, then 
to correctly define Gromov-Witten invariant out of the moduli space $\overline{\cal M}_{g,k}(M, d[{\bf P}^1])$, we 
have to take acount of the Euler class of the obstruction bundle as described above.

So we want to consider the integrals where the class $\Phi$ appears in the integrand. Formally, we have

{\bf Definition 4.1:} $K_{(k,g,d,\Theta)} := \int_{\overline{\cal M}_{g,k}({\bf P}^n,d)}\Theta\wedge \Phi$, 
where $\Theta$ is a polynomial in Chern classes of certain equivariant vector bundles over 
$\overline{\cal M}_{g,k}({\bf P}^n,d)$.

For example, let us consider the case of mukai flop. It is well known that the normal bundle of the embeded ${\bf P}^n$ 
is actaully naturally isomorphic to its cotangent bundle because of the existence of holomorphic 2-forms.

{\bf Definition 4.2:} $K_{(3,0,d,ev^*(\alpha) \wedge ev^*(\beta) \wedge ev^*(\gamma))}:= \int_{\overline{\cal M}_{0,3}({\bf P}^n,d)} ev^*(\alpha)\wedge ev^*(\beta)\wedge ev^*(\gamma) \wedge \Phi$.
where $\alpha$,$\beta$,$\gamma$ are any cohomogy classes of ${\bf P} ^n$ with appropriate degrees, i.e. 
\begin {equation}
\deg (\alpha)+ \deg(\beta)+ \deg(\gamma)+ \deg(\Phi) = \dim {\cal M}_{0,3}({\bf P}^n,d).
\end {equation}

Note that the invariant defined above includes all qauntum correction coming from exceptional effective curve in the case of mukai flop.

About these invariants, we have the following vanishing theorem

{\bf Theorem 4.3:} The invariants $K_{(k,g,d, \Theta)}$ all vanish regardless of the flexibility of $\Theta$.

\subsection{Bott's residue formula and normal bundle contibution}

In this subsection, we introduce the technique we use to compute the invariants as defined in previous subsections.
 The basic ideal is to consider torus action and use the Bott's residue formula to reduce the integral to fixed points
 loci of the action.
 
 Starting from \cite{K}, a lot of work has been done towards localization techniques applied to the computation of
 Gromov-Witten invariants and verification of mirror symmetry predictions. In the most general case, one has to consider
 localization of virtual classes as done in \cite{GP, LLY}. In \cite{CKYZ}, the authors developed effective ways to compute
 similar invariants involving Euler classes of obstruction bundles. But they mainly treat direct sums of line bundles.
 For our computation, the machinery introduced by \cite{K} suffices. Here we will follow the presentation in \cite{K}
 closely. To keep notation simple, we will only consider integration formula in genus zero case. The proof of vanishing 
 of the invariants in higher genus case will be almost identical. We will point out the slight difference later.
 
 Before proving theorem $4.3$, we first want to introduce {\bf Bott's residue formula:}
 
 Let $X$ be a compact complex projective manifold (orbifold allowed) and $E$ a holomorphic vector bundle (or orbibundle)
 over $X$. Suppose $T:=(C^*)^{n+1}$ a complex torus acts on $(X,E)$. Denote the fixed points loci by $X^T$ and its
 connected components by $X^{\gamma}$. Since the irreducible representations of torus are dimensional one, over 
 $X^{\gamma}$ the bundle $E$ splits into direct sum of line bundles $E^{\gamma, \lambda}$ twisted by character 
 $\lambda : T \longrightarrow C^*$, $\lambda\in T^\vee = Z \oplus Z \oplus\cdots\oplus Z$. The normal bundle of 
 $X^{\lambda}$ (denoted by $N^{\lambda}$) also splits into sum of line bundles $N^{\gamma, \lambda }$ over characters
 $\lambda\in T^\vee\setminus \{0 \}$.
 
 By splitting principle, we suppose the Chern classes of bundle $E$ are given by homogeneous symmetric polynomials
 in degree 2 generators $e_i$'s as follows:
 \begin{equation}
 \sum_{k \geq 0}c_k(E) = \Pi_i(1+e_i), \,\,\,\,\,\, e_i \in H^2(X, Q).
 \end{equation}
 
 Analogously, we add generators $e_i^{\gamma, \lambda}$ and $n_i^{\gamma, \lambda}$ to $H^2(X^{\gamma}, Q)$.
 
 Let $ P$ be a homogeneous symmetric polynomial. Then the {\bf Bott's residue formula} reads:
 \begin{equation}
   \int_X P(e_i) = \sum_{\gamma}\int_{X^\gamma}\frac{P(e_i^{\gamma,\lambda} + \lambda)}{\Pi(n_i^{\gamma,\lambda})}.
 \end{equation}
 
 The right hand side of the above formula is considered as rational function in $\lambda$'s. It turns out to have 
 homogeneous degree 0 ( actually a constant independent of choice of $\lambda$'s). The numerator of r.h.s. is actually
 the equivariant extension of the pullback of class $P(e_i)$ to $X^\gamma$. The denominator is the equivariant Euler
 class of normal bundle of $X^\gamma$.
 
 Now, we want to calculate the fixed points in the moduli space of stable maps in order to apply the Bott's residue 
 formula.
 
 Let $T= (C^*)^{n+1}$ acts diagonally on ${\bf P}^n$ with generic weights $-\lambda_1, -\lambda_2,\cdots,-\lambda_{n+1}$.
 The fixed points are projectivization of coordinate lines of $C^{n+1}$, denoted by $p_i$. And the only invariant curves 
 are lines connecting the fixed points labeled by $l_{ij} = l_{ji}$, where $i\not= j$.
 
 The action of $T$ on ${\bf P}^n$ induces an action of $T$ on the moduli space of stable maps $\overline{\cal M}_{g,k}({\bf P}^n,d)$
 by moving the image of the map. Let $(C,f,x_1,\cdots, x_k)$ be a fixed point in the stable map space. Then the geometric
 image of the map is fixed. So we have
 \begin{enumerate}
 \item The contracted components, the marked points, the ramification points, the nodes all are mapped to the fixed 
 points $p_i$'s in ${\bf P}^n$.
 \item A non-contracted component is map onto one of the lines $l_{ij}$'s, ramifying over the two fixed points(end 
 points of the line ), thus is forced to be rational and completely determined by its degree. 
\end{enumerate}

We associate with each fixed 
 map a marked graph $\Gamma$ as follows. The vertices of the graph $v \in Vert(\Gamma)$ correspond to the connected 
 components $C_v$ of $f^{-1}(p_1,p_2,\cdots, p_{n+1})$. Here the component can be either a point or union of 
 irreducible components of the curve $C$. The edges $\alpha\in Edge(\Gamma)$ correspond to non-contracted
 component of $C^\alpha$ of genus $0$ mapping onto the $l_{ij}$'s. There are also tails on the vertices coming from the
 marked points. We also mark the graph by the following labels:
 
 \begin{enumerate}
 \item Label the vertices numbers $f_v$ from $1$ to $n+1$ defined by $f(C_v):= p_{f_v}$. Also label a vertex by $g_v$
 (the genus of the $1$-dimensional part of $C_v$, for a point the genus is $0$) and a set $S_v \subset \{1,2,\cdots,k \}$
 the indices of the marked points.
 \item Label the edges by the mapping degree $d_{\alpha}\in N$
 \end{enumerate}
  
 The claim is that the connected components of  $\overline{\cal M}_{g,k}({\bf P}^n,d)^T$ are isomorphic to 
 $\Pi_{v\in Vert(\Gamma)} \overline{\bf M}_{g_v,val(v)}/\mbox{\bf Aut}(\Gamma)$ and can be identified as
  equivalent classes of connected graphs $\Gamma$ with labeling satisfying the following
  conditions:

  \begin{enumerate}
  \item[(1)] For $\alpha\in Edge(\Gamma)$ connecting vertices $u,v \in  Vert(\Gamma)$, then $f_u \not= f_v$,
  \item[(2)] $1 - \chi (\Gamma) + \sum_{v \in Vert(\Gamma)} g_v = g$,
  \item[(3)] $\sum_{\alpha\in Edge(\Gamma)} d_{\alpha} = d$,
  \item[(4)] $\cup_{v \in Vert(\Gamma)}S_v = \{1,2,\cdots,k \}$.
  \end{enumerate}

   From now on we only consider the integration formula for genus $0$ case. We first want
   to give some notations:
   \begin{enumerate}
   \item[(1)] For a graph, we define an incident pair of vertex and edge
   $(v, \alpha)$ to be a flag $F = (v,\alpha)$ and denote by $w_F$ the
   expression $\frac{\lambda_{f_v}-\lambda_{f_u}}{d_\alpha}$ where $u\not=v$
   is the other vertex of the edge $\alpha$.
   \item[(2)] Recall that $\overline{\cal M}_{0,k}$ is the Deligne-Mumford space
   of marked stable curves. For each marking $i$, there is a line bundle
   $L_i\longrightarrow \overline{\cal M}_{0,k}$ with fiber $T^*_{C,p_i}$ over the
   moduli point $C$. Define $\psi_i := c_1(L_i)$.
   \end{enumerate}

Now we describe the normal bundle of the fixed points components. For an
equivariant bundle $E$, denote by $[E]$ its class in the corresponding equivariant
K-group. Also we denote $\overline{\cal M}_{0,k}({\bf P}^n,d)$ by
$\overline{\cal M}$ for simplicity and often denote a bundle by its geometric
fiber at a point $(C,f)$.

To keep notation simple, we ignore the marked points as in \cite{K} and explain
the difference along the way.

The class of normal bundle for a component $\overline{\cal M}^{\gamma}$ having
graph type $\Gamma$ is

\begin{equation}
  [N_{\overline{\cal M}^\gamma}]  =  [T_{\overline{\cal M}}] - [T_{\overline{\cal M}^\gamma}]
\end{equation}
$$
  [T_{\overline{\cal M}}]  =  [H^0(C,f^*(T{\bf P}^n))] + \sum_{y \in C^\alpha\cap C^\beta} [T_y(C^\alpha)\otimes T_y(C^\beta)] \nonumber 
$$
\begin{equation}
    +  \sum_{y \in C^\alpha\cap C^\beta : \alpha\not= \beta}([T_y(C^\alpha)] + [T_y(C^\beta)]) - \sum_\alpha[H^0(C^\alpha,TC^\alpha)]
\end{equation}

The first summand corresponds to infinitesmal deformation of the map $f$ from $C$. The second summand corresponds to
smoothing of nodes. And the third comes from deformation of the curve $C$ fixing the singular points. If there is a 
marked point $x$ on $C^\alpha$, it should also be fixed and in the third summand there would be an additional term
$ \sum_\alpha[H^0(T_x(C^\alpha))]$. (Same remark applies to the formula below).

\begin{eqnarray}
[T_{\overline{\cal M}^\gamma}] & = & \sum_{y \in C^\alpha\cap C^\beta: \alpha\not=\beta:\alpha,\beta \not\in Edge(\Gamma)}
[T_y(C^\alpha)\otimes T_y(C^\beta)]\nonumber \\
 & + & \sum_{y \in C^\alpha\cap C^\beta: \alpha\not=\beta:\alpha \not\in Edge(\Gamma)}
 [T_y(C^\alpha)] - \sum_{\alpha \not\in Edge(\Gamma)}[H^0(C^\alpha, T C^\alpha)].
\end{eqnarray}  
where the first term corresponds to smoothing of nodes which are intersection of two contracted components. 
The second term and the third come from deformation of the components preserving singular points.

So we have the following formula 
\begin{equation}
  [N_{\overline{\cal M}^\gamma}] = [H^0(C,f^*(T{\bf P}^n))]+ [N_{\overline{\cal M}^\gamma}^{abs}]
\end{equation}
where
$$
[N_{\overline{\cal M}^\gamma}^{abs}] := \sum_{y \in C^\alpha\cap C^\beta: \alpha\not=\beta:\alpha,\beta \in Edge(\Gamma)}
[T_y(C^\alpha)\otimes T_y(C^\beta)]
$$
$$
+ \sum_{y \in C^\alpha\cap C^\beta: \alpha\in Edge(\Gamma),\beta \not\in Edge(\Gamma)}
[T_y(C^\alpha)\otimes T_y(C^\beta)]
$$
\begin{equation}
 + \sum_{y \in C^\alpha\cap C^\beta: \alpha\not=\beta:\alpha\in Edge(\Gamma)}[T_y(C^\alpha)]
 - \sum_{\alpha\in  Edge(\Gamma)}[H^0(C^\alpha,T C^\alpha)].
\end{equation}

In the formula for $ [N_{\overline{\cal M}^\gamma}^{abs}]$ above the first and third summand are trivial bundles twisted 
with characters of the torus. The term $[H^0(C,f^*(T{\bf P}^n))]$ and the classes from the bundle $E$ restricted to 
$X^\gamma$ in our application later have same nature. When we take the Chern classes of these summand, we just get weights 
of torus action on the fibers of these bundles (expressed in terms of $\lambda_i$'s), hence can be pulled out of the 
integral. In the second summand, the tangent space of the non-contracted component at $y$ is fixed but twisted, while 
the tangent space of the contracted component at $y$ is moving without twisting. Taking equivariant Chern class we get 
a sum of certain tangential weight and the $\psi$ class over suitable space of pointed stable curves. This reduce the 
integral on the right hand side of Bott's formula integral to integral of $\psi$ classes over space of pointed curves 
for which the answer has been conjectured by Witten and verified by Kontsevich rigorously. Thus we have a contribution 
(as rational function in $\lambda$'s) from each of the admissible graphs. The invariant is given by a graph summation 
collecting all these contributions:

\begin{eqnarray}
  & &\prod_{\alpha \in Edge(\Gamma );v_1, v_2:vertices \ of \ \alpha}
      (\frac{(-1)^{d_\alpha}(\frac{d_\alpha}{\lambda _{V_1}-\lambda _{V_2}})^{2 d_\alpha}}{(d_\alpha !)^2}) \nonumber \\
  & \times & \prod_{\alpha \in {Edge(\Gamma )}}\prod_{k \not= f_{v_1},k \not= f_{v_2}}\prod_{a,b \geq
     0:a+b=d_{\alpha }}\frac{1}{\frac{a}{d_\alpha }\lambda _{f_{v_1}}+\frac{b}{d_\alpha }\lambda _{f_{v_2}}-\lambda _k}\nonumber \\
  & \times & \prod_{v \in Vert(\Gamma)}\{(\sum_{flags:F=(v,\alpha )}w_F^{-1})^{ val (v)-3}\times 
  \prod_{ flags:F=(v,\alpha )}w_F^{-1}\nonumber  \\
  &\times & \prod_{j \not= f_v}(\lambda _{f_v}-\lambda _j)^{val (v)-1}\}.
\end{eqnarray}  

Here the valence of a vertex includes the counts of the number of
tails. The detailed calculation of the weights can be found in \cite{K} which 
we refer the interested readers to.

\subsection{Proof of the vanishing theorem.}

 In this subsection, we prove the vanishing theorem stated in Section
4.2. We show the calculation for the specific example defined in Definition 
{\bf 4.2 } with n=2. The proof for the general cases is almost identical. We will briefly explain the difference at the end of the proof.

{\bf Proof of Theorem 4.3:} First of all, note that the invariant is given by
\begin{eqnarray}
& & \sum _{\Gamma }\frac{1}{|\mbox{Aut}(\Gamma )|}\times
(\mbox{contribution from }{ev}^*(\alpha)\wedge {ev}^*(\beta) \wedge {ev}^*(\gamma) )\nonumber  \\
& \times & (\mbox{contribution from} \Phi )\times (\mbox{formula} (25)).
\end{eqnarray}

Here the contribution of the second and the third terms are just a product of the weights of induced torus action 
on the corresponding vector bundles. We show that the contribution from the Euler class $\Phi $ of the obstruction
bundle restricted to the fixed point component is zero and thus conclude.

To deal with nodal curves, we need the following normalization sequence.
First let us consider the simple case where $\bf C=C_\alpha \cup C_\beta $. There is an exact sequence of maps 
of sheaves (of the holomorphic functions):
\begin{equation}
 0 \longrightarrow {\cal O}_C \longrightarrow {\cal O}_{C_\alpha}\oplus {\cal O}_{C_\beta}\longrightarrow
 {\cal O}_{ C_\alpha\cap C_\beta}\longrightarrow 0. 
\end{equation}

Here all the maps except the last one are obtained from inclusions.\ And the
last one maps $(f_1,f_2)$ to $f_1-f_2$.

In general we have the normalization sequence resolving all the
nodes of ${\bf C}$ which are forced by a graph type $\Gamma $
\begin{eqnarray}
 0 \longrightarrow {\cal O}_C & \longrightarrow & \left(\oplus _{v \in Vert(\Gamma )}{\cal O}_{C_v}\right) 
 \oplus \left( \oplus _{\alpha \in Edge(\Gamma )}{\cal O}_{C_\alpha }\right) \nonumber \\
 & \longrightarrow & \oplus _{F \in Flag(\Gamma )}{\cal O}_{x_F}\longrightarrow  0,
 \end{eqnarray}
where $x_F=C_v \cap C_\alpha$ for a flag $(v,\alpha )$, and the last map
sends $(g|_{C_V},h|_{C_\alpha})$ to $g-h$ on the intersection point.

Twist the above sequence by $f^*T^*{\bf P}^2$ and take cohomology to get
\begin{eqnarray}
0 & \longrightarrow & H^0(C,f^*T^*{\bf P}^{2})  \nonumber  \\
& \longrightarrow & \left( \oplus _{v \in Vert(\Gamma )}H^0(C_v,f^*T^*{\bf P}^{2})\right) \oplus 
\left( \oplus _{\alpha \in Edge(\Gamma )}H^0(C_\alpha,f^*T^*{\bf P}^{2})\right)\nonumber \\
 & \longrightarrow & \oplus _{F \in Flag(\Gamma )}T_{f(x_F)}^*{\bf P}^{2}\longrightarrow H^1(C,f^*T^*{\bf P}^{2})
 \nonumber \\
 & \longrightarrow & \left( \oplus _{v \in Vert(\Gamma )}H^1(C_v,f^*T^*{\bf P}^{2})\right) \oplus 
 \left( \oplus _{\alpha \in Edge(\Gamma )}H^1(C_\alpha,f^*T^*{\bf P}^{2})\right)\nonumber \\ 
 & \longrightarrow & 0.
 \end{eqnarray} 
The first term in the third line follows since $x_F$ is a point, which is why the last term in the last line is 0. 
Note that $f^*T^*{\bf P}^2|_{C_v}$ is trivial since $C_v$ is mapped to a point , hence $H^0(C_v,f^*T^*{\bf P}^2)
=T^*{\bf P}^2|_{P_{f(v)}}$ and $H^1(C_v,f^*T^*{\bf P}^2)=H^1(C_v,{\cal O})\otimes f^*T^*{\bf P}^2$. Since we are 
considering genus zero case, $H^1(C_v,f^*T^*{\bf P}^2)$ is also zero( In general, it can be expressed in terms 
of the first Chern class $C_1$ of Hodge bundle over space of pointed curves).  Because of the concavity
of $T^*{\bf P}^2$, $H^0(C_\alpha,f^*T^*{\bf P}^2)$ is zero. And by looking
at the maps in the first line of $(29)$, we see $H^0(C,f^*T^*{\bf P}^2)$ is also zero. So we have 
\begin{equation}
 [H^1(C,f^*T^*{\bf P}^2)]=[\prod _{v \in Vert(\Gamma )}T_{p_{f(v)}}^*{\bf P}^2]+[\prod _{a \in Edge(\Gamma )}H^1
 (C_\alpha,f^*T^*{\bf P}^2)].
\end{equation} 
The contribution from the l.h.s is the product of those of the two terms on
the r.h.s.

Since a non-contracted component is rigid, $[H^1(C_\alpha,f^*T^*{\bf P}^2)]$ is a trivial bundle
when restricted to fixed point components. To compute the weights, we
consider the following description of the cotangent bundle of ${\bf P}^2$
by an exact sequence of bundles over ${\bf P}^2={\bf P}(V)$ where $V$ is a complex vector space of dimension 3. 
First we have
\begin{equation}
0 \longrightarrow {\cal O}(-1)\longrightarrow V \longrightarrow Q \longrightarrow 0,
\end{equation}
where $V$ represents the trivial bundle with vector space $V$ as fiber and  ${\cal O}(-1)$ is the universal bundle. 
Tensoring with ${\cal O}(1)$, we have
\begin{equation}
0 \longrightarrow {\cal O}\longrightarrow {\cal O}(1)\otimes V \longrightarrow {\cal O}(1)\otimes 
Q \longrightarrow 0,
\end{equation}
where ${\cal O}(1)\otimes Q ={\bf T P}^2$. Dualizing we have
\begin{equation}
0 \longrightarrow T^*{\bf P}^2 \longrightarrow 
{\cal O}(-1)\otimes V^*\longrightarrow {\cal O}\longrightarrow 0.
\end{equation}

Pulling back by $f$ over $C_\alpha$ and taking cohomology, we have
\begin{eqnarray}
0 & \longrightarrow & H^0(C_\alpha,f^*T^*{\bf P}^2)\longrightarrow H^0(C_\alpha,{\cal O}(-d)\otimes V^*)\nonumber \\
 & \longrightarrow & H^0(C_\alpha,{\cal O})\longrightarrow H^1(C_\alpha,f^*T^*{\bf P}^2) \nonumber \\
& \longrightarrow & H^1(C_\alpha,{\cal O}(-d)\otimes V^*)\longrightarrow H^1(C_\alpha,{\cal O})\longrightarrow 0.
\end{eqnarray}
Note that $C_\alpha $ is rational. $H^0(C_\alpha,{\cal O}(-d)\otimes V^*)$ and $H^1(C_\alpha,{\cal O})$ are both $0$. 
So we have
\begin{equation}
 0 \longrightarrow H^0(C_\alpha,{\cal O})\longrightarrow H^1(C_\alpha,f^*T^*{\bf P}^2)\longrightarrow
 H^1(C_\alpha,{\cal O}(-d)\otimes V^*)\longrightarrow 0.
 \end{equation}
So the contribution of $[H^1(C_\alpha,f^*T^*{\bf P}^2)]$ is given by a product of weights on 
$H^1(C_\alpha$, ${\cal O}(-d)\otimes V^*)$ and weight on $H^0(C_\alpha,{\cal O})$.
Obviously the weight on $H^0(C_\alpha,{\cal O})$  is zero. So the contribution of $[H^1(C_\alpha, f^*T^*{\bf P}^2)]$ 
is zero for each $\alpha \in Edge(\Gamma )$.

From $(30)$, we see that the total contribution of the Euler class $\Phi $ is zero. Thus we conclude our proof 
of the genus zero case.

In the general case of higher genus, formula $(26)$ needs to be modified. $\lambda$ classes (coming from 
the deformation of the complex structures on ${\bf C})$ and hence Hodge integrals will appear in the computation 
of the normal bundle contribution and the details can be found in \cite{GP}. But the point is that the
contribution of Euler class $\Phi $ is still zero, since there is a 0 weight coming from $H^1(C_\alpha,f^*T^*{\bf P}^2)$
 for each non-contracted component (necessarily rational as explained earlier). So Theorem $ 4.3$ still holds.

{\bf Theorem 4.4:} Suppose that non-singular projective manifolds $X$ and $X'$ of complex dimension $2n$
 are connected by a sequence of Mukai flops. Then $X$ and $X' $ have isomorphic Ruan cohomologies.
 
 {\bf Proof:} By theorem 4.3, we have that all Gromov-Witten invariants appearing in the right hand side of $(1)$ 
 vanish. Therefore, we have that for $X, X'$ their quantum corrections all vanish. Thus their quantum cohomology
 are the same as their ordinary Chow ring. By theorem 3.2, we know that $X, X'$ have isomorphic Ruan cohomology. This 
 proves the theorem.
 
 {\bf Corollary 4.5:} For Mukai flops, cohomological minimal model conjecture holds.
 
 Finally, we present a well known proposition to point out that local existence of a holomorphic symplectic  2-form implies
 natural isomorphism of the normal bundle and the cotangent bundle for a embedded ${\bf P}^n$.
 
 {\bf Proposition 4.6: (see \cite{Mukai})} Suppose that ${\bf P}^n$ is embedded in a smooth variety $X$ with a neighborhood $N$ admitting
 a holomorphic symplectic 2-form $\omega$, then we have the following
 \begin{enumerate}
 \item $ \mbox{codim}_X{\bf P}^n \geq n$.
 \item In case $\mbox{codim}_X{\bf P}^n = n$, there is a natural isomorphism $T^*{\bf P}^n = N_{X|{\bf P}^n}$.
 \end{enumerate} 
 
 {\bf Proof:} Since $H^{2,0}({\bf P}^n) = 0$, $\omega \mid _{{\bf P}^n} = 0$. Thus $T_p{\bf P}^n 
 \subset (T_p{\bf P}^n)^\bot$ for any point $p \in {\bf P}^n$, where $T_p{\bf P}^n \subset T_pX$ is considered as a
 subspace of $T_pX$. Hence codim$_X{\bf P}^n = \dim (T_p{\bf P}^n)^\bot \geq \dim T_p{\bf P}^n = n $. In case equality
 holds, $T_p{\bf P}^n = (T_p{\bf P}^n)^\bot$.
 
 $\omega \mid _{T_pX}$ is nondegenerate, so there is an isomorphism $\phi : T_p X = (T_p X)^*$. Thus we have 
 $T_p{\bf P}^n =(T_p{\bf P}^n)^\bot = \mbox{Ann}(T_p{\bf P}^n) = N^*_{X \setminus{\bf P}^n}$, where the second isomorphism 
 is via the map $\phi$.

Email address: stsjxhu@zsu.edu.cn

\end{document}